\documentclass[12pt,reqno]{amsart}
\usepackage{amssymb}
\usepackage{amscd}
\usepackage{hyperref}

\newtheorem{theorem}{Theorem}[section]

\newtheorem{proposition}[theorem]{Proposition}

\numberwithin{equation}{section}
\begin{document}
\title{The spectrum of basic Dirac operators}
\author[K.~Richardson]{Ken Richardson}
\address{Department of Mathematics \\
Texas Christian University \\
Fort Worth, Texas 76129, USA}
\email[K.~Richardson]{k.richardson@tcu.edu}

\begin{abstract}
In this note, we discuss Riemannian foliations, which are smooth foliations
that have a transverse geometric structure. We explain a known generalization of
Dirac-type operators to transverse operators called ``basic Dirac
operators'' on Riemannian foliations, which require
the additional structure of what is called a bundle-like metric. We explain the result in 
\cite{HabRi} that the spectrum of such an operator
is independent of the choice
of bundle-like metric, provided that the transverse geometric structure is
fixed. We discuss consequences, which include defining a new version of the
exterior derivative and de Rham cohomology that are nicely adapted to this
transverse geometric setting.
\end{abstract}

\maketitle

\section{Introduction}

The content here concerns some work in \cite{HabRi} and also briefly
mentions work in \cite{PaRi}, \cite{LRi1}, and \cite{BKR}; it
also provides applications not given in the these references.

\subsection{Smooth foliations and basic forms}

Let $(M,\mathcal{F})$ be a smooth, closed manifold of dimension $n$ endowed
with a foliation $\mathcal{F}$ given by an integrable subbundle $L\subset TM$
of rank $p$. The set $\mathcal{F}$ is a partition of $M$ into immersed
submanifolds (\emph{leaves}) such that the transition functions for the
local product neighborhoods (foliation charts) are smooth. The subbundle $L=T%
\mathcal{F}$ is the tangent bundle to the foliation; at each $p\in M$, $T_{p}%
\mathcal{F}=L_{p}$ is the tangent space to the leaf through $p$.

\vspace{0in}Many researchers have studied basic forms and basic cohomology,
especially in the particular cases of Riemannian foliations with bundle-like
metrics, to be discussed later (see \cite{Al}, \cite{KT1}, \cite{To}). Basic
forms are differential forms on $M$ that locally depend only on the
transverse variables in the foliation charts --- that is, forms $\alpha $
satisfying $X\lrcorner \alpha =X\lrcorner d\alpha =0$ for all $X\in \Gamma
(L)$; the symbol \textquotedblleft $\lrcorner $\textquotedblright\ stands
for interior product. Let $\Omega \left( M,\mathcal{F}\right) \subset \Omega
\left( M\right) $ denote the space of basic forms. These differential forms
are preserved by the exterior derivative and are used to define basic
cohomology groups, which can be infinite-dimensional but are always
finite-dimensional in the case of Riemannian foliations. We define the \emph{%
basic cohomology group }$H^{k}\left( M,\mathcal{F}\right) $ by%
\begin{equation*}
H^{k}\left( M,\mathcal{F}\right) =\frac{\ker d_{k}}{\mathrm{Im}d_{k-1}}
\end{equation*}%
with%
\begin{equation*}
d_{k}=d:\Omega ^{k}\left( M,\mathcal{F}\right) \rightarrow \Omega
^{k+1}\left( M,\mathcal{F}\right) .
\end{equation*}%
We comment that the point of using basic forms is an effort to find a form
of de Rham cohomology on a singular, possibly non-Hausdorff space, that
space being the set of leaves of the foliation. To gain the smooth
structure, we loose a bit of information about the leaf space. The basic
cohomology can be infinite-dimensional, and it can be relatively trivial. We
may also define basic cohomology with values in a foliated vector bundle; by
doing this we gain more topological information about the leaf space.

Basic cohomology does not necessarily satisfy Poincar\'{e} duality, even if
the foliation is transversally oriented. If there are additional
restrictions, it may satisfy duality, for example if the manifold admits a
metric for which the leaves are locally equidistant and are minimal
submanifolds (ie a taut Riemannian foliation). We emphasize that basic
cohomology is a smooth foliation invariant and does not depend on the choice
of metric or any transverse or leafwise structure.

\subsection{Riemannian foliations and bundle-like metrics%
}

We assume throughout the paper that the foliation is \emph{Riemannian}; this
means that there is a metric on the local space of leaves --- a
holonomy-invariant transverse metric $g_{Q}$ on the normal bundle $%
Q=TM\diagup L$; this means that the transverse Lie derivative $\mathcal{L}%
_{X}g_{Q}$ is zero for all leafwise vector fields $X\in \Gamma (L)$. This
metric is a substitute for a metric on the singular space of leaves. This
condition is characterized by the existence of a unique metric and
torsion-free connection $\nabla $ on $Q$ \cite{Mo}, \cite{Re}, \cite{To}. We
can then associate to $\nabla $ the transversal curvature data, in
particular the transversal Ricci curvature $\mathrm{Ric}^{\nabla }$ and
transversal scalar curvature $\mathrm{Scal}^{\nabla }$.

We often assume that the manifold is endowed with the additional structure
of a \emph{bundle-like metric} \cite{Re}, i.e. the metric $g$ on $M$ induces
the metric on $Q\simeq L^{\perp }$. Every Riemannian foliation admits
bundle-like metrics that are compatible with a given $\left( M,\mathcal{F}%
,g_{Q}\right) $ structure. There are many choices, since one may freely
choose the metric along the leaves and also the transverse subbundle $N%
\mathcal{F}$. We note that a bundle-like metric on a smooth foliation is
exactly a metric on the manifold such that the leaves of the foliation are
locally equidistant.

There are topological restrictions to the existence of bundle-like metrics
(and thus Riemannian foliations). Important examples of requirements for the
existence of a Riemannian foliations include

\begin{itemize}
\item certain characteristic classes must vanish (see \cite{KT2})

\item leaf closures must partition the manifold (see \cite{Mo})

\item the basic cohomology must be finite-dimensional (see \cite{KT1}, \cite%
{To}, \cite{PaRi})

\item for any metric on the manifold, the orthogonal projection 
\begin{equation*}
P:L^{2}\left( \Omega \left( M\right) \right) \rightarrow L^{2}\left( \Omega
\left( M,\mathcal{F}\right) \right)
\end{equation*}%
must map the subspace of smooth forms onto the subspace of smooth basic
forms (\cite{PaRi}).
\end{itemize}

Riemannian foliations were introduced by B. Reinhart in 1959 (\cite{Re}).
Good references for Riemannian foliations and bundle-like metrics include
the books and papers of B. Reinhart, F. W. Kamber, Ph. Tondeur, P. Molino,
for example.

\subsection{The basic Laplacian}

Many researchers have studied basic forms and the basic Laplacian on
Riemannian foliations with bundle-like metrics (see \cite{Al}, \cite{KT1}, 
\cite{To}). The basic Laplacian $\Delta _{b}$ for a given bundle-like metric
is a version of the Laplace operator that preserves the basic forms and that
is essentially self-adjoint on the $L^{2}$-closure of the space of basic
forms. We define the basic Laplacian $\Delta _{b}$ by%
\begin{equation*}
\Delta _{b}=d\delta _{b}+\delta _{b}d:\Omega \left( M,\mathcal{F}\right)
\rightarrow \Omega \left( M,\mathcal{F}\right) ,
\end{equation*}%
where $\delta _{b}$ is the $L^{2}$-adjoint of the restriction of $d$ to
basic forms: $\delta _{b}=P\delta $ is the ordinary adjoint of $d$ followed
by the orthogonal projection onto the space of basic forms.

The operator $\Delta _{b}$ and its spectrum depend on the choice of the
bundle-like metric and provide invariants of that metric. See \cite{JuRi}, 
\cite{LRi1}, \cite{LRi2}, \cite{PaRi}, \cite{Ri1}, \cite{Ri2} for results.
One may think of this operator as the Laplacian on the space of leaves. This
operator is the appropriate one for physical intuition. For example, the
Laplacian is used in the heat equation, which determines the evolution of
the temperature distribution over a manifold as a function of time. If we
assume that the leaves of the foliation are perfect conductors of heat, then
the basic Laplacian is the appropriate operator that allows one to solve the
heat distribution problem in this situation.

It turns out that the basic Laplacian is the restriction to basic forms of a
second order elliptic operator on all forms, and this operator is not
necessarily symmetric (\cite{PaRi}). Only in special cases is this operator
the same as the ordinary Laplacian.

The basic Laplacian $\Delta _{b}$ is also not the same as the formal
Laplacian defined on the local quotient manifolds of the foliation charts
(or on a transversal). This transversal Laplacian is in general not
symmetric on the space of basic forms, but it does preserve $\Omega \left( M,%
\mathcal{F}\right) $.

The basic heat flow asymptotics are more complicated than that of the
standard heat kernel, but there is a fair amount known (see \cite{PaRi}, 
\cite{Ri1}, \cite{Ri2}).

\subsection{The basic adjoint of the exterior derivative and mean curvature}

We assume $\left( M,\mathcal{F},g_{M}\right) $ is a Riemannian foliation
with bundle-like metric compatible with the Riemannian structure $\left( M,%
\mathcal{F},g_{Q}\right) $. For later use, we define the mean curvature
one-form $\kappa $ and discuss the operator $\kappa _{b}\lrcorner $. Let 
\begin{equation*}
H=\sum_{i=1}^{p}\pi \left( \nabla _{f_{i}}^{M}f_{i}\right) ,
\end{equation*}%
where $\pi :TM\rightarrow N\mathcal{F}$ is the bundle projection and $\left(
f_{i}\right) _{1\leq i\leq p}$ is a local orthonormal frame of $T\mathcal{F}$%
. This is the mean curvature vector field, and its dual one-form is $\kappa
=H^{\flat }$. Let $\kappa _{b}=P\kappa $ be the (smooth) basic projection of
this mean curvature one-form. Let $\kappa _{b}\lrcorner $ denote the
(pointwise) adjoint of the operator $\kappa _{b}\wedge $. Clearly, $\kappa
_{b}\lrcorner $ depends on the choice of bundle-like metric $g_{M}$, not
simply on the transverse metric $g_{Q}$.

It turns out that $\kappa _{b}$ is a closed form whose cohomology class in $%
H^{1}\left( M,\mathcal{F}\right) $ is independent of the choice of
bundle-like metric (see \cite{Al}).

Recall the following expression for $\delta _{b}$, the $L^{2}$-adjoint of $d$
restricted to the space of basic forms of a particular degree (see \cite{To}%
, \cite{PaRi}):%
\begin{eqnarray*}
\delta _{b} &=&P\delta \\
&=&\pm \overline{\ast }d\overline{\ast }+\kappa _{b}\lrcorner \\
&=&\delta _{T}+\kappa _{b}\lrcorner ,
\end{eqnarray*}%
where

\begin{itemize}
\item $\delta _{T}$ is the formal adjoint (with respect to $g_{Q}$) of the
exterior derivative on the transverse local quotients.

\item the pointwise transversal Hodge star operator $\overline{\ast }$ is
defined on all $k$-forms $\gamma $ by%
\begin{equation*}
\overline{\ast }\gamma =\left( -1\right) ^{p\left( q-k\right) }\ast \left(
\gamma \wedge \chi _{\mathcal{F}}\right) ,
\end{equation*}%
with $\chi _{\mathcal{F}}$ being the leafwise volume form, the
characteristic form of the foliation and $\ast $ being the ordinary Hodge
star operator. Note that $\overline{\ast }^{2}=\left( -1\right) ^{k\left(
q-k\right) }$ on $k$-forms.

\item The sign $\pm $ above only depends on dimensions and the degree of the
basic form.
\end{itemize}

\subsection{Twisted duality for basic cohomology}

Even for transversally oriented Riemannian foliations, Poincar\'{e} duality
does not necessarily hold for basic cohomology.

However, note that $d-\kappa _{b}\wedge $ is also a differential which
defines a cohomology of basic forms. That is, since $d\left( \kappa
_{b}\right) =0$, it follows from the Leibniz rule that $\left( d-\kappa
_{b}\wedge \right) ^{2}=0$ as an operator on forms, and it maps basic forms
to basic forms. This differential also has the property that%
\begin{equation*}
\delta _{b}\overline{\ast }\alpha =\left( -1\right) ^{k+1}\overline{\ast }%
\left( d-\kappa _{b}\wedge \right) \alpha 
\end{equation*}%
on every basic $k$-form $\alpha $. As a result, the transversal Hodge star
operator implements an isomorphism between different kinds of basic
cohomology groups (see \cite{KTduality}, \cite{To}, and \cite{PaRi}):%
\begin{equation*}
H_{d}^{\ast }\left( M,\mathcal{F}\right) \cong H_{d-\kappa _{b}\wedge
}^{q-\ast }\left( M,\mathcal{F}\right) .
\end{equation*}%
This is called \emph{twisted Poincar\'{e} duality}.

\subsection{Ordinary Dirac operators and examples}

See a reference such as \cite{Roe} for the well-known details for this
section. The ordinary Dirac operator in Euclidean space is given by 
\begin{equation*}
D=\sum e_{k}\cdot \frac{\partial }{\partial x_{k}},
\end{equation*}%
where the operators $e_{k}\cdot $ are multiplication by matrices satisfying
the relation 
\begin{equation*}
e_{k}\cdot e_{j}\cdot ~+~e_{j}\cdot e_{k}\cdot ~=-2\delta _{kj}\mathbf{1}.
\end{equation*}%
In three dimensions, the matrices can be chosen to be the Pauli spin matrices%
\begin{equation*}
e_{1}\cdot ~=\left( 
\begin{array}{cc}
0 & -1 \\ 
1 & 0%
\end{array}%
\right) ,~e_{2}\cdot ~=\left( 
\begin{array}{cc}
0 & i \\ 
i & 0%
\end{array}%
\right) ,~e_{3}\cdot ~=\left( 
\begin{array}{cc}
i & 0 \\ 
0 & -i%
\end{array}%
\right) .
\end{equation*}

These relations are the same as the relations of the complex Clifford
algebra $\mathbb{C}\mathrm{l}(V)$ associated to a vector space $V$. In
general what is needed to define an ordinary Dirac operator on a Riemannian
manifold $M$ is a vector bundle $E\rightarrow M$ that is a bundle of $%
\mathbb{C}\mathrm{l}(TM)$ Clifford modules with compatible connection $%
\nabla ^{E}$. The\emph{\ Dirac operator} $D$ is the composition of the maps 
\begin{equation*}
\Gamma \left( E\right) \overset{\nabla ^{E}}{\longrightarrow }\Gamma \left(
T^{\ast }M\otimes E\right) \overset{\cong }{\longrightarrow }\Gamma \left(
TM\otimes E\right) \overset{\mathrm{Cliff}}{\longrightarrow }\Gamma \left(
E\right) ,
\end{equation*}%
where the last map is Clifford multiplication, denoted by \textquotedblleft $%
\cdot $\textquotedblright . We may write 
\begin{equation*}
D=\sum e_{i}\cdot \nabla _{e_{i}}^{E}~
\end{equation*}%
acting on $\Gamma \left( E\right) $, where $\left( e_{i}\right) $ is a local
orthonormal frame of $TM$. Computations show that $D$ is elliptic,
essentially self-adjoint, and thus has discrete spectrum.

Examples of the Dirac operator include:

\begin{itemize}
\item \textquotedblleft The\textquotedblright\ spin$^{c}$ Dirac operator.
Here $E$ is a spinor bundle, which at each point $p\in M$ is an irreducible
representation space for $\mathbb{C}\mathrm{l}(T_{p}M)$. 

\item The de Rham operator%
\begin{equation*}
d+\delta :\Omega ^{\mathrm{even}}\left( M\right) \rightarrow \Omega ^{%
\mathrm{odd}}\left( M\right) .
\end{equation*}

\item The signature operator%
\begin{equation*}
d+\delta :\Omega ^{+}\left( M\right) \rightarrow \Omega ^{-}\left( M\right) ,
\end{equation*}%
where the $\pm $ refer to self-dual and anti-self-dual forms. There is an
operator of the form $\bigstar =i^{k\left( k-1\right) +n}\ast $ acting on
complex-valued $k$-forms on a $2n$-dimensional manifold, and $\bigstar ^{2}=%
\mathbf{1}$. The space $\Omega ^{+}$ of self-dual forms is the $+1$
eigenspace of $\bigstar $, and the space $\Omega ^{-}$ of antiself-dual
forms is the $-1$-eigenspace of $\bigstar $. It turns out that $d+\delta $
anticommutes with $\bigstar $ and thus maps $\Omega ^{\pm }$ to $\Omega
^{\mp }$.

\item The Dolbeault operator%
\begin{equation*}
\partial +\overline{\partial }:\Omega ^{0,\mathrm{even}}\left( M\right)
\rightarrow \Omega ^{0,\mathrm{odd}}\left( M\right) .
\end{equation*}
\end{itemize}

Each one of these operators has an associated Laplacian $D^{2}$ and
associated harmonic forms, and the index of each of these operators ( index$%
\left( D\right) =\dim \ker D-\dim \ker D^{\ast }$ ) is an important
topological invariant. For example, if $D=d+\delta $ is the de Rham
operator, we have%
\begin{equation*}
\ker \left( d+\delta \right) =\mathcal{H,}
\end{equation*}%
the space of harmonic forms, which by Hodge theory can be used to represent
the different cohomology classes. Thus, 
\begin{eqnarray*}
\text{index}\left. \left( d+\delta \right) \right\vert _{\Omega ^{\mathrm{%
even}}} &=&\dim \ker \left. \left( d+\delta \right) \right\vert _{\Omega ^{%
\mathrm{even}}}-\dim \ker \left. \left( d+\delta \right) \right\vert
_{\Omega ^{\mathrm{odd}}} \\
&=&\chi \left( M\right) ,
\end{eqnarray*}%
the Euler characteristic of $M$.

\subsection{The basic Dirac operator and statement of the main theorem}

We now discuss the construction of the basic Dirac operator (see \cite{DGKY}%
, \cite{GlK}, \cite{PrRi}, \cite{BKR}), a construction which requires a choice of
bundle-like metric. Let $(M,\mathcal{F})$ be a Riemannian manifold endowed
with a Riemannian foliation. Let $E\rightarrow M$ be a foliated vector
bundle (see \cite{KT2}) that is a bundle of $\mathbb{C}\mathrm{l}(Q)$
Clifford modules with compatible connection $\nabla ^{E}$. The \emph{%
transversal Dirac operator} $D_{\mathrm{tr}}$ is the composition of the maps 
\begin{equation*}
\Gamma \left( E\right) \overset{\left( \nabla ^{E}\right) ^{\mathrm{tr}}}{%
\longrightarrow }\Gamma \left( Q^{\ast }\otimes E\right) \overset{\cong }{%
\longrightarrow }\Gamma \left( Q\otimes E\right) \overset{\mathrm{Cliff}}{%
\longrightarrow }\Gamma \left( E\right) ,
\end{equation*}%
where the last map denotes Clifford multiplication, denoted by
\textquotedblleft $\cdot $\textquotedblright , and the operator $\left(
\nabla ^{E}\right) ^{\mathrm{tr}}$ is the projection of $\nabla ^{E}$. The
transversal Dirac operator fixes the basic sections $\Gamma _{b}(E)\subset
\Gamma (E)$ (i.e. $\Gamma _{b}(E)=\{s\in \Gamma (E):\nabla _{X}^{E}s=0$ for
all $X\in \Gamma (L)\}$) but is not symmetric on this subspace. By modifying 
$D_{\mathrm{tr}}$ by a bundle map, we obtain a symmetric and essentially
self-adjoint operator $D_{b}$ on $\Gamma _{b}(E)$. Let $\kappa _{b}$ be the $%
L^{2}$-orthogonal projection of $\kappa $ onto the space of basic forms as
explained above, and let $\kappa _{b}^{\sharp }$ be the corresponding vector
field. We now define 
\begin{eqnarray*}
D_{\mathrm{tr}}~s &=&\sum_{i=1}^{q}e_{i}\cdot \nabla _{e_{i}}^{E}s~, \\
D_{b}s &=&\frac{1}{2}(D_{\mathrm{tr}}+D_{\mathrm{tr}}^{\ast
})s=\sum_{i=1}^{q}e_{i}\cdot \nabla _{e_{i}}^{E}s-\frac{1}{2}\kappa
_{b}^{\sharp }\cdot s~,
\end{eqnarray*}%
where $\{e_{i}\}_{i=1,\cdots ,q}$ is a local orthonormal frame of $Q$. A
direct computation shows that $D_{b}$ preserves the basic sections, is
transversally elliptic, and thus has discrete spectrum (\cite{GlK}, \cite%
{DGKY}, \cite{HabRi}).

An example of the basic Dirac operator is as follows. Using the bundle $%
\wedge ^{\ast }Q$ as the Clifford bundle with Clifford action $e\cdot
~=e^{\ast }\wedge -e^{\ast }\lrcorner $ in analogy to the ordinary de Rham
operator, we have%
\begin{eqnarray*}
D_{\mathrm{tr}} &=&d+\delta _{T}=d+\delta _{b}-\kappa _{b}\lrcorner :\Omega
^{\mathrm{even}}\left( M,\mathcal{F}\right) \rightarrow \Omega ^{\mathrm{odd}%
}\left( M,\mathcal{F}\right)  \\
D_{b} &=&\frac{1}{2}(D_{\mathrm{tr}}+D_{\mathrm{tr}}^{\ast })s=d+\delta
_{b}-\kappa _{b}\lrcorner -\frac{1}{2}\kappa _{b}^{\sharp }\cdot  \\
&=&d+\delta _{b}-\kappa _{b}\lrcorner -\frac{1}{2}\left( \kappa _{b}\wedge
-\kappa _{b}\lrcorner \right)  \\
&=&d+\delta _{b}-\frac{1}{2}\kappa _{b}\lrcorner -\frac{1}{2}\kappa
_{b}\wedge .
\end{eqnarray*}%
One might have instead guessed that $d+\delta _{b}$ is the basic de Rham
operator in analogy to the ordinary de Rham operator, for this operator is
essentially self-adjoint, and the associated basic Laplacian yields basic
Hodge theory that can be used to compute the basic cohomology.

We study the invariance of the spectrum of the basic Dirac operator with
respect to a change of bundle-like metric; that means when one modifies the
metric on $M$ in any way that leaves the transverse metric on the normal
bundle intact (this includes modifying the subbundle $N\mathcal{F}\subset TM$%
, as one must do in order to make the mean curvature basic, for example). In 
\cite{HabRi}, we prove

\begin{theorem}
\label{inv} Let $(M,\mathcal{F})$ be a compact Riemannian manifold endowed
with a Riemannian foliation and basic Clifford bundle $E\to M$. The spectrum
of the basic Dirac operator is the same for every possible choice of
bundle-like metric that is associated to the transverse metric on the
quotient bundle $Q$.
\end{theorem}

We emphasize that the basic Dirac operator $D_{b}$ depends on the choice of
bundle-like metric, not merely on the Clifford structure and Riemannian
foliation structure, since both projections $T^{\ast }M\rightarrow Q^{\ast }$
and $P_{b}$ as well as $\kappa _{b}$ depend on the leafwise metric. 

\section{Proof of the main theorem}

The proof of the main theorem is contained in \cite{HabRi}. The idea of
proof is as follows. One can show that every different choice of bundle-like
metric changes the $L^{2}$-inner product by multiplication by a specific
smooth, positive basic function. This changes the basic Dirac operator by a
zeroth order operator that is Clifford multiplication by an exact basic
one-form. This new operator is conjugate to the original one, and thus the
spectrum of the operator is independent of the metric choice.

\section{Consequences of the main theorem}

D. Dominguez showed that every Riemannian foliation admits a bundle-like
metric for which the mean curvature form is basic \cite{Do}. Further, the
bundle-like metric may be chosen so that the mean curvature is
basic-harmonic (in the new metric); see \cite{MarMinooRuh} and \cite{Ma}.
Therefore, in calculating or estimating the eigenvalues of the basic Dirac
operator, one may choose the bundle-like metric so that the mean curvature
is basic-harmonic. Immediately we may obtain stronger inequalities for
eigenvalue estimates. We give one example below.

In \cite{Ju}, S. D. Jung showed that the eigenvalues $\lambda $ of the basic
Dirac operator on spin foliations (where the normal bundle carries a spin
structure as a foliated vector bundle) satisfy 
\begin{equation*}
\lambda ^{2}\geq \frac{q}{4(q-1)}\mathrm{inf}_{M}(\mathrm{Scal}^{\nabla
}+\left\vert \kappa \right\vert ^{2}),
\end{equation*}%
under the assumption that the mean curvature form $\kappa $ is basic and
basic-harmonic. In \cite{Hab}, the author obtained another Friedrich-type
estimate \cite{Fr} for the eigenvalues of the basic Dirac operator for
bundle-like metrics with basic-harmonic mean curvature. Since by Theorem \ref%
{inv} the spectrum does not change, we may improve both results to deduce
that for \textbf{any} bundle-like metric the eigenvalues of the basic Dirac
operator satisfy 
\begin{equation}
\lambda ^{2}\geq \frac{q}{4(q-1)}\mathrm{inf}_{M}(\mathrm{Scal}^{\nabla }).
\label{eq:esti}
\end{equation}%
This estimate is of interest only for positive transversal scalar curvature.
Moreover, in \cite{HabRi} we show

\begin{proposition}
\label{pro:esti} Let $(M,\mathcal{F})$ be a compact Riemannian manifold
endowed with a spin foliation with basic mean curvature $\kappa.$ Then, we
have the estimate 
\begin{equation*}
\lambda^2\geq \frac{q}{4(q-1)}\mathrm{inf}_M(\mathrm{Scal}_M-\mathrm{Scal}%
_L+|A|_Q^2+|T|_L^2).  \label{eq:estima}
\end{equation*}
If $\mathcal{F}$ is a Riemannian flow (i.e. $p=1$), then 
\begin{equation}
\lambda^2\geq \frac{q}{4(q-1)}\mathrm{inf}_M(\mathrm{Scal}%
_M+|A|_Q^2+|\kappa|^2).  \label{eq:estmflot}
\end{equation}
If the limiting case is attained, the foliation is minimal and we have a
transversal Killing spinor.
\end{proposition}

Here $A$ and $T$ denote the O'Neill tensors \cite{B,O} of the foliation.

More applications can be found in the paper \cite{HabRi}.

\section{Modified differentials, Laplacians, and basic cohomology}

\vspace{0in}From the above, the basic de Rham operator is 
\begin{eqnarray*}
D_{b} &=&d+\delta _{b}-\frac{1}{2}\kappa _{b}\lrcorner -\frac{1}{2}\kappa
_{b}\wedge  \\
&=&\widetilde{d}+\widetilde{\delta }
\end{eqnarray*}%
acting on basic forms, where%
\begin{equation*}
\widetilde{d}=d-\frac{1}{2}\kappa _{b}\wedge ,\widetilde{\delta }=\delta
_{b}-\frac{1}{2}\kappa _{b}\lrcorner .
\end{equation*}%
The operators $\widetilde{d}$ and $\widetilde{\delta }$ have interesting
properties:

\textbullet They are differentials: $\widetilde{d}^{2}=0$, $\widetilde{%
\delta }^{2}=0$.

\textbullet $\widetilde{\delta }\overline{\ast }=\pm \overline{\ast }%
\widetilde{d}$.

From the first property we see that we can define cohomology using these
differentials. From the second property, we know that we may define a basic
signature operator. This was not known and not possible previously with
ordinary basic cohomology and the operator $d+\delta _{b}$, because that
operator would not map self-$\overline{\ast }$-dual basic forms to anti-self-%
$\overline{\ast }$-dual forms. The eigenvalues of this operator depend only
on the Riemannian foliation structure.

The new basic cohomology defined using $\widetilde{d}$ will satisfy Poincar%
\'{e} duality, and the isomorphism is implemented using $\overline{\ast }$.
Also, even though the differential depends on the choice of the bundle-like
metric, the dimensions of the resulting cohomology groups are independent of
that choice. The consequences of these results will be explained in a
forthcoming paper.

\end{document}